\documentclass[12pt,a4paper]{amsart}

\usepackage[english]{babel}
\usepackage[T1]{fontenc}
\usepackage[utf8]{inputenc}

\usepackage{amsmath, amsthm, amscd, amsfonts}
\usepackage{amssymb, tikz}
\usepackage{mathtools}

\usepackage{hyperref}
\usepackage[capitalise]{cleveref}

\usepackage{setspace}
\newcommand{\MPB}{\mathbb{P}}
\newcommand{\MQB}{\mathbb{Q}}
\newcommand{\MRB}{\mathbb{R}}

\DeclareMathOperator{\cf}{cf}

\newcommand{\Ult}{\mathrm{Ult}}

\newcommand{\crit}{\mathrm{crit}}

\newcommand{\Add}{\mathrm{Add}}
\newcommand{\Sacks}{\mathrm{Sacks}}

\newcommand{\Col}{\mathrm{Col}}

\newcommand{\lusim}[1]{\smash{\underset{\raisebox{1.2pt}[0cm][0cm]{$\sim$}}
{{#1}}}}

\newtheorem{theorem}{Theorem}[section]
\newtheorem{lemma}[theorem]{Lemma}

\theoremstyle{definition}
\newtheorem{definition}[theorem]{Definition}

\theoremstyle{remark}
\newtheorem{remark}[theorem]{Remark}
\newtheorem{claim}[theorem]{Claim}

\title[No universal graphs at uncountable regular cardinals]{No universal graphs at uncountable regular cardinals}

\author[M. Golshani]{Mohammad  Golshani}
\thanks{{\color{red}This is the preliminary version of the paper which only sketches the proof. The details will appear in the next version
of the paper.}}

\thanks{%
 The author's research has been supported by a grant from
  IPM (No. 1401030417). He thanks Omer Ben-Neria for some useful comments.}

\address{
School of Mathematics, Institute for Research in Fundamental Sciences (IPM), P.O. Box:
19395-5746, Tehran-Iran. }
\email{golshani.m@gmail.com}
\urladdr{http://math.ipm.ac.ir/~golshani/}

\subjclass[]{}


\begin{document}
\begin{abstract}
Assuming the existence of a strong cardinal, we find a model of ZFC in which for each uncountable regular cardinal $\lambda,$
there is no universal graph of size $\lambda$.
\end{abstract}
 \maketitle

\section{Introduction}
The existence of universal objects in an infinite cardinal is of great interest in both set theory and model theory. In this paper we concentrate on universal graphs. It is well-known that there is a countable universal graph and that if $\lambda$ is an uncountable  cardinal
with $2^{<\lambda}=\lambda$, then there exists a universal graph of size $\lambda$.
 On the other hand, by a result of Shelah (see \cite{kojman}), if $\lambda$  is a regular cardinal with $2^\lambda=\lambda^+$, then after forcing with $\Add(\lambda, \lambda^{++})$ for adding $\lambda^{++}$-many Cohen subsets to $\lambda,$ there are no universal graphs of size $\lambda^+.$ Friedman and Thompson \cite{friedman-thompson} obtained further results in this direction, in particular, they showed it is consistent, modulo the existence of strong cardinals, that there are no universal graphs at the successor of a singular strong limit cardinal of countable cofinality.
In this paper we continue this line of research further and prove the following global result:
\begin{theorem}
\label{maintheorem}
Assume GCH holds and $\kappa$ is a $(\kappa+3)$-strong cardinal. Then there exists a generic extension  $W$ of the universe in which
$\kappa$ remains inaccessible and for each uncountable regular cardinal $\lambda < \kappa,$ there is no universal graph of size $\lambda,$
in particular the rank initial segment of $W$ at $\kappa$ is a model of
ZFC in which for each uncountable regular cardinal $\lambda,$ there is no universal graph of size $\lambda.$
\end{theorem}
The theorem in particular answers Question 6.1 from \cite{friedman-thompson}.
The rest of the paper is devoted to the proof of the above theorem. In Section \ref{s0} we present some preliminaries. In Section \ref{s1}
 we construct the model $W$, and then in Section \ref{s2}
we show that $W$ is as required.

\section{Some preliminaries}
\label{s0}

We will start by stating some general non-existing results for universal graphs.
\begin{lemma} (Shelah, see \cite{kojman})
\label{shtm}
Suppose $\lambda$ is a regular cardinal and $2^\lambda=\lambda^+.$
Then in the generic extension by $\Add(\lambda, \mu)$, for $\cf(\mu) \geq \lambda^{++},$
there are no universal graphs of size $\lambda^+.$
\end{lemma}
Friedman and Thompson \cite{friedman-thompson} isolated the main properties used in the proof of the above lemma and
they have given an abstract lemma extending the above lemma, we will need the following special case.
\begin{lemma} (see \cite{friedman-thompson})
\label{sdftm}
Suppose $\lambda$ is a regular cardinal and $2^\lambda=\lambda^+.$
Then in the generic extension by $\Sacks(\lambda, \mu)$, for $\cf(\mu) \geq \lambda^{++},$
there are no universal graphs of size $\lambda^+.$
\end{lemma}

\section{Building the model W}
\label{s1}
In this section we construct the generic extension $W$ which will be the desired model as requested in Theorem \ref{maintheorem}. We construct the model $W$ in essentially two steps. The first step is a reverse Easton iteration of suitable forcing notions of length $\kappa+1$.
 The second step is a variant of Radin forcing using the existence of suitable guiding generics.

Assume GCH holds, $\kappa$ is a $(\kappa+3)$-strong cardinal and let $j:V \to M$ witness $\kappa$ is $(\kappa+3)$-strong. Let also
$i: V \to N$ be the ultrapower embedding by $U=\{X \subseteq \kappa: \kappa \in j(X)    \}$. Also factor $j$ through $i$ to get
$k: N \to M$ defined by $k([f]_U)=j(f)(\kappa)$. Note that $\crit(k)=\kappa^{++}_N < \kappa^{++}_M=\kappa^{++}$.

 Force over $V$ with the reverse Easton iteration
\[
\langle \langle  \MPB_\lambda: \lambda \leq \kappa+1       \rangle, \langle \lusim{\MQB}_\lambda: \lambda \leq \kappa   \rangle\rangle
\]
 where we force with the trivial forcing notion except  $\lambda\leq \kappa$ is an inaccessible cardinal in which case

$$\Vdash_{\MPB_\lambda}\text{``} \lusim{\MQB}_\lambda=\lusim{\Sacks}(\lambda, \lambda^{++}) \times \lusim{\Add}(\lambda^{+}, \lambda^{+3}) \times \lusim{\Add}(\lambda^{++}, \lambda^{+4}).$$
Let $\langle  \langle  G_\lambda: \lambda \leq \kappa+1 \rangle, \langle    H_\lambda: \lambda \leq \kappa       \rangle\rangle$
be generic for the forcing, so each $G_\lambda$ is $\MPB_\lambda$-generic over $V$ and $H_\lambda$
is $\MQB_\lambda=\lusim{\MQB}_\lambda[G_{\lambda}]$-generic over $V[G_\lambda]$. Let us also write, for each inaccessible $\lambda \leq \kappa,$
$H_\lambda=H^{0}_\lambda \times H^{1}_\lambda \times H^{2}_\lambda$, which corresponds to
$\MQB_\lambda=\Sacks(\lambda, \lambda^{++}) \times \Add(\lambda^{+}, \lambda^{+3}) \times \Add(\lambda^{++}, \lambda^{+4}).$

By arguments similar to \cite{friedman-honzik} and \cite{cummings}, the following hold in $V[G_{\kappa+1}]$:
\begin{enumerate}
\item $\kappa$ is a $(\kappa+2)$-strong cardinal,

\item  $j: V \to M$ extends to some $j^*: V[G_{\kappa+1}] \to M[G_{\kappa+1} \ast H]$
with $^{\kappa}M[G_{\kappa+1} \ast H] \subseteq M[G_{\kappa+1} \ast H], V_{\kappa+3}[G_{\kappa+1}] \subseteq M[G_{\kappa+1} \ast H],$
and $j^*$ is generated by a $(\kappa, \kappa^{+3})$-extender,

\item $i: V \to N$ extends to $i^*: V[G_{\kappa+1}] \to N[i^*(G_{\kappa+1})]$ and $i^*$ is the ultrapower of $V[G_{\kappa+1}]$
by $U^*=\{X \subseteq \kappa: \kappa \in j^*(X)   \}$

\item For each inaccessible cardinal $\lambda \leq \kappa,$ we have
\[
2^\lambda=\lambda^{++}+2^{\lambda^+}=\lambda^{+3}+2^{\lambda^{++}}=\lambda^{+4},
\]

\item $M[G_{\kappa+1} \ast H] \models``2^\kappa=\kappa^{++}+2^{\kappa^+}=\kappa^{+3}+2^{\kappa^{++}}=\kappa^{+4}$'',

\item There is an $F\in V[G_{\kappa+1}]$ which is a generic filter over $N[i^*(G_{\lambda+1})]$, by the forcing notion
$$(\Col(\kappa^{+4}, < i(\kappa)) \times \Add(\kappa^{+3}, i(\kappa)) \times \Add(\kappa^{+4}, i(\kappa)^+))_{N[i^*(G_{\kappa+1})]}$$

\end{enumerate}
Note that
$F \in M[G_{\kappa+1} \ast H]$, as it can be coded by an element of $V_{\kappa+2}$, also,
if $k^*:  N[i^*(G_{\kappa+1})] \to M[G_{\kappa+1} \ast H]$ is the induced elementary embedding, defined by
$k^*([f]_{U^*})=j^*(f)(\kappa),$
then $\crit(k^*)=\kappa^{++}_{N} < \kappa^{++}_{M}=\kappa^{++}$
and $F$ can be transferred along $k^*$, in the sense that $\langle k^*{}''[F] \rangle$, the filter generated by $k^*{}''[F]$, is
$(\Col(\kappa^{+4}, < i(\kappa)) \times \Add(\kappa^{+3}, i(\kappa)) \times \Add(\kappa^{+4}, i(\kappa)^+))_{M[G_{\kappa+1} \ast H]}$-generic over
$M[G_{\kappa+1} \ast H].$

 For notational simplicity let us denote the models $V[G_{\kappa+1}],\ M[G_{\kappa+1} \ast H]$
 and $N[i^*(G_{\kappa+1})]$ by $V^*, M^*$ and $N^*$ respectively.

 Work in $V^*$. Let
 \[
 \mathcal{U}=\langle \mathcal{U}(\alpha, \beta): \alpha \leq \kappa, \beta < o^{\mathcal U}(\alpha) \rangle
 \]
 be a coherent sequence of measures of length $\ell^{\mathcal U}=\kappa+1$
 and $o^{\mathcal U}(\kappa)=\kappa^+$.
 For each such $\alpha, \beta$ let
 \[
 i^*_{\alpha, \beta}: V^* \to M^*_{\alpha, \beta} \simeq \Ult(V^*, \mathcal{U}(\alpha, \beta))
 \]
 be the induced embedding.
  Let also
 \[
\mathcal{K}= \langle   K_{\alpha, \beta}: \alpha \leq \kappa, \beta < o^{\mathcal U}(\alpha)         \rangle
 \]
 be a sequence of filters such that:
 \begin{enumerate}
 \item[(7)] $K_{\alpha, \beta}$ is $\bigg(\Col(\kappa^{+4}, < i_{\alpha, \beta}(\kappa)) \times \Add(\kappa^{+3}, i_{\alpha, \beta}(\kappa)) \times  \Add(\kappa^{+4}, i_{\alpha, \beta}(\kappa)^+)\bigg)_{M^*_{\alpha, \beta}}$-generic over $M^*_{\alpha, \beta}$,

 \item[(8)] the sequence is coherent in the sense that
 \[
 \langle K_{\alpha, \tau}: \tau < \beta  \rangle = [\bar\alpha \mapsto \langle K_{\bar\alpha, \tau}: \tau < \beta \rangle]_{\mathcal{U}(\alpha, \beta)}.
 \]
 \end{enumerate}
 For each $\alpha \leq \kappa$ set $\mathcal{F}(\alpha)=\bigcap_{\beta < o^{\mathcal U}(\alpha)} \mathcal{U}(\alpha, \beta)$
 if $o^{\mathcal U}(\alpha) >0,$ and set $\mathcal{F}(\alpha)=\{\emptyset\}$ otherwise.
 We now use $\mathcal{U}$ and $\mathcal{K}$
 to define a variant of Radin forcing.
 \begin{definition}
 Let $\MRB=\MRB_{\mathcal{U}, \mathcal{K}}$
 be the forcing notion consisting of conditions of the form
 \[
 p=\langle \alpha^p_{-1}, (\alpha^p_0, f^p_0, A^p_n, F^p_0), \cdots, (\alpha^p_{n^p}, f^p_{n^p}, A^p_{n^p}, F^p_{n^p})                  \rangle
 \]
 where
 \begin{enumerate}
 \item $n^p\geq 0,$

 \item $\alpha^p_{-1} < \alpha^p_0 < \cdots < \alpha^p_{n^p}=\kappa$,


 \item  $f^p_i \in \Col((\alpha^p_{i-1})^{+4}, < \alpha^p_i) \times \Add((\alpha^p_{i-1})^{+3}, \alpha^p_i) \times \Add((\alpha^p_{i-1})^{+4}, (\alpha^p_i)^+),$

 \item $A^p_i \in \mathcal{F}(\alpha^p_i)$,

 \item $F^p_i$ is a function with domain $A^p_i$ such that:
 \begin{enumerate}
 \item[(A)] if $o^{\mathcal{U}}(\alpha^p_i) > 0$, then:
 \begin{enumerate}
 \item for each $\theta \in A^p_i, \theta > \alpha^p_{i-1}$,

 \item $f^p_i \in \Col((\alpha^p_{i-1})^{+4}, < \theta) \times \Add((\alpha^p_{i-1})^{+3}, \theta) \times \Add((\alpha^p_{i-1})^{+4}, \theta^+)$, for each $\theta \in A^p_i,$

 \item for all $\theta \in A^p_i$, $F^p_i(\theta) \in \Col(\theta^{+4}, < \alpha^p_i) \times \Add(\theta^{+3}, \alpha^p_i) \times \Add(\theta^{+4}, (\alpha^p_i)^+),$

 \item $[F^p_i]_{\mathcal{U}(\alpha^p_i, \beta)} \in K_{\alpha^p_i, \beta}$, for all $\beta < \o^{\mathcal{U}}(\alpha)$.
 \end{enumerate}
 \item[(B)] if $o^{\mathcal{U}}(\alpha^p_i) = 0$, then $F^p_i(\emptyset)=f^p_i.$
 \end{enumerate}
 \end{enumerate}
 The order relations $\leq$ and $\leq^*$ (the Prikry order) on $\MRB$ are defined in the natural way, see for example \cite{cummings} or \cite{golshani}.
 \end{definition}
Let $\MRB=\MRB_{\mathcal{U}, \mathcal{K}}$ and let $K$ be $\MRB$-generic over $V^*$.
Let $C$ be the Radin club added by $K$ and let $\langle \lambda_i: i<\kappa \rangle$
be an increasing enumeration of $C$. For each $i<\kappa$
let $K_i=K_i^C \times K_i^3 \times K_i^4 $ be $\Col(\lambda_i^{+4}, <\lambda_{i+1}) \times \Add(\lambda_i^{+3}, \lambda_{i+1}) \times \Add(\lambda_i^{+4}, \lambda^+_{i+1})$-generic filer over $V^*$ added by $K$.
In the next lemma we collect the main properties of the above forcing notion.
 \begin{lemma}
 \label{rforcingprop}
 \begin{enumerate}
 \item (the chain condition) $\MRB$ satisfies the $\kappa^+$-Knaster,

 \item (Prikry property) The forcing notion $(\MRB, \leq, \leq^*)$ satisfies the Prikry property.

 \item (the factorization lemma): Suppose $p \in \MRB$ is of length $n> 0$ and let $0 \leq m < n$. Then
 \begin{enumerate}
\item  there exists an isomorphism
  $$\MRB / p \simeq \MRB / p^{\leq m} \times \MRB / p^{>m},$$
  where
 $p^{\leq m}=p \restriction m+1= \langle \alpha^p_{-1}, (\alpha^p_0, f^p_0, A^p_n, F^p_0), \cdots, (\alpha^p_{m}, f^p_{m}, A^p_{m}, F^p_{m})                  \rangle$
 and
 $p^{>m}= \langle \alpha^p_m, (\alpha^p_{m+1}, f^p_{m+1}, A^p_{m+1}, F^p_{m+1}), \cdots, (\alpha^p_{n^p}, f^p_{n^p}, A^p_{n^p}, F^p_{n^p})                  \rangle$,

 \item $(\MRB/ p^{\leq m}, \leq)$
 satisfies the $(\alpha^p_m)^+$-c.c. and $(\MRB / p^{>m}, \leq^*)$ is $(\alpha^p_m)^{+3}$-closed.

 \end{enumerate}
 \item $\kappa$ remains an inaccessible cardinal in $V^*[K]$.

 \item Let $C$ and $\langle \lambda_i: i<\kappa \rangle$ be as above. Then
 \begin{enumerate}
 \item for each $i<\kappa$, we have
 \[
 2^{\lambda_i}=\lambda_i^{++}+2^{\lambda_i^+}=\lambda_i^{+3}+2^{\lambda_i^{++}}=
 \lambda_i^{+4}+2^{\lambda_i^{+3}}=\lambda_i^{+5}=\lambda_{i+1}+2^{\lambda_i^{+4}}=\lambda_{i+1}^+.
 \]
 \item If we force with $\Col(\aleph_0, < \lambda_0) \times \Add(\aleph_0, \lambda_0^+) \times \MRB$,
 then in the generic extension, $\kappa$ is the least inaccessible cardinal and
 for all $\lambda < \kappa, 2^\lambda=\lambda^{++}$.
 \end{enumerate}

 \item $V^*[K]=V^*[\langle K_i: i<\kappa \rangle       \rangle]$.
 \end{enumerate}
 \end{lemma}
 Now suppose that $L_0 \times L_1 \times K$ is $\Col(\aleph_0, < \lambda_0) \times \Add(\aleph_0, \lambda_0^+) \times \MRB$-generic over
 $V^*$
 and set
 \[
 W=V^*[L_0 \times L_1 \times K] = V[G_{\kappa+1}][L_0 \times L_1 \times K].
 \]
 This completes our construction of the model $W$. Let us note that
   by Lemma \ref{rforcingprop},
\[
\text{Card}^W \cap [\aleph_0, \kappa)=\{\aleph_0\} \cup \{\lambda_i, \lambda^+_i, \lambda^{++}_i, \lambda^{+3}_i, \lambda^{+4}_i: i<\kappa                         \}.
\]

In the next section we show that the model $W$ satisfies the conclusion of Theorem \ref{maintheorem}.
\section{No uncountable universal graphs in W}
\label{s2}
In this section we complete the proof of Theorem \ref{maintheorem}, by showing that in $W$,
there is no universal graph of size $\lambda$, where $\lambda$ is an uncountable regular cardinal less than
$\kappa$.

The proof is divided into several cases:

{\underline{\bf Case 1. $\lambda=\aleph_1$:}}

We have
$W=V^*[K, L_0][L_1]$, where $V^*[K, L_0]$ satisfies CH and $W$ is an extension of $V^*[K, L_0]$ by $L_1$
which is $\Add(\aleph_0, \lambda_0^+)$-generic over $V^*[K, L_0]$. Thus by Lemma \ref{shtm}, there are no
universal graphs of size $\aleph_1$.

{\underline{\bf Case 2. $\lambda=\lambda_{i+1}^{+n}$ for some $i < \kappa$ and some $n < 3$:}}

Assume on the contrary that $T$ is a universal graph of size $\lambda$ in $W$.  Now we can write $W$ as
\[
W=V^*[\langle  K_j: j \leq i     \rangle][\langle K_j: i < j <\kappa      \rangle][L_0 \times L_1]
\]
and that by Lemma \ref{rforcingprop}(2), $\mathcal{P}^W(\lambda) = \mathcal{P}^{V^*[\langle  K_j: j \leq i    \rangle][L_0 \times L_1]}(\lambda)$.
It immediately follow that $T \in V^*[\langle  K_j: j \leq i      \rangle][L_0 \times L_1].$

Now let $i_*$ be a limit ordinal, with possibly $i_*=0$, such that $i=i_*+k$, for some $k < \omega.$
So we have
\[
V^*[\langle  K_j: j \leq i      \rangle][L_0 \times L_1]=V^*[\langle  K_j: j < i_*      \rangle][L_0 \times L_1][\prod_{\ell \leq k}K_{i_*+\ell}].
\]
Now
$V^*=V[G_{\lambda_{i+1}}][X]$, where $X$ is generic over $V[G_{\lambda_{i+1}}]$ for a forcing notion which is $\lambda_{i+1}^*$-closed
and $\lambda_{i+1}^*$ is the least inaccessible above $\lambda_{i+1}$. It immediately follow
that $\mathcal{U} \restriction \lambda_{i_*}+1 \in V[G_{\lambda_{i+1}}]$ and hence by Lemma \ref{rforcingprop}(2),
we may assume that $T$ belongs to
$$V[G_{\lambda_{i+1}}][\langle  K_j: j \leq i      \rangle][H_{\lambda_{i+1}}][L_0 \times L_1]=V[G_{\lambda_{i+1}}][\langle  K_j: j \leq i      \rangle][L_0 \times L_1][\prod_{\ell < 3, \ell \neq n }H^\ell_{\lambda_{i+1}}][ H^n_{\lambda_{i+1}}].$$
Now we have $V[G_{\lambda_{i+1}}][\langle  K_j: j \leq i      \rangle][L_0 \times L_1][\prod_{\ell < 3, \ell \neq n }H^\ell_{\lambda_{i+1}}]$
satisfies $2^{\lambda}=\lambda^{+}$
and $ V[G_{\lambda_{i+1}}][\langle  K_j: j \leq i      \rangle][H_{\lambda_{i+1}}][L_0 \times L_1]$
is a generic extension of it, by either $\Sacks(\lambda, \lambda^{++})$ if $\lambda=\lambda_{i+1}$
or $\Add(\lambda, \lambda^{++})$ if $\lambda \in \{\lambda_{i+1}^+, \lambda_{i+1}^{++}  \}$. This contradicts either Lemma \ref{shtm} or Lemma
\ref{sdftm}.

{\underline{\bf Case 3. $\lambda=\lambda_{i}^{+n}$ for some $i < \kappa$ and some $2 < n < 5$:}}

 The argument is essentially the same as above. Suppose towards a contradiction that $T$ is a universal graph of size $\lambda$ in $W$. As in Case 2, we will get that $T \in V[G_{\lambda_{i}+1}][\langle  K_j: j < i      \rangle][K_{i}][L_0 \times L_1]$ and this model is an extension of $V[G_{\lambda_{i}+1}][\langle  K_j: j < i      \rangle][K^C_{i} \times K^{\ell}_{i}][L_0 \times L_1]$, where $\ell \in \{3, 4\}, \ell \neq n$. This later model satisfies GCH at $\lambda$,
and as in Case 2 we get a contradiction to Lemma \ref{shtm}.

{\underline{\bf Case 4. $\lambda=\lambda_{i}^{+n}$ for some limit ordinal $i < \kappa$ and some $0< n < 3$:}}

Arguing as above, if $T \in W$ is universal graph of size $\lambda$,
then $T \in V[G_{\lambda_i}][\langle K_j: j<i \rangle][H^0_{\lambda_i} \times H^1_{\lambda_i} \times H^2_{\lambda_i}][L_0 \times L_1]$
and this model extends $V[G_{\lambda_i}][\langle K_j: j<i \rangle][H^0_{\lambda_i} \times H^\ell_{\lambda_i}][L_0 \times L_1]$,
where $\ell \in \{1, 2\}, \ell \neq n$, by $\Add(\lambda, \lambda^{++})$, so again we get a contradiction
using Lemma \ref{shtm}.

{\underline{\bf Case 5. $\lambda=\lambda^+_i$ for some limit ordinal $i< \kappa$:}}

This is the hardest part of the proof. We follow very closely the argument
given in \cite{friedman-thompson}. Suppose towards a contradiction that $T$ is a universal graph on  $\lambda$. As before, we can conclude that
$T \in V[G_{\lambda_i}][\langle K_j: j<i \rangle][L_0 \times L_1][H^0_{\lambda_i} \times H^1_{\lambda_i}]$.
Let us identify $H^0_{\lambda_i}$ with the $\kappa^{++}$-sequences
$\langle  A_\sigma: \sigma < \kappa^{++}          \rangle$ of mutually Sacks-generic subsets of $\lambda_i$.

For notational simplicity, set
\[
\widetilde{V}=V[G_{\lambda_i}][\langle K_j: j<i \rangle][L_0 \times L_1][H^1_{\lambda_i}],
\]
and for each $\sigma \leq \lambda_i^{++}$ set
\[
\widetilde{V}_\sigma=V[G_{\lambda_i}][\langle K_j: j<i \rangle][L_0 \times L_1][H^1_{\lambda_i}][\langle A_\alpha: \alpha < \sigma    \rangle].
\]
Thus $\widetilde{V}=\widetilde{V}_0$ and $\widetilde{V}_{\lambda_i^{++}}=V[G_{\lambda_i}][\langle K_j: j<i \rangle][L_0 \times L_1][H^0_{\lambda_i} \times H^1_{\lambda_i}].$ Also for each $\sigma < \lambda_i^{++}$,
the model $\widetilde{V}_\sigma$ satisfies $2^{\lambda_i}=\lambda_i^+.$
The proof of the next claim is essentially the same as in \cite{mohammadpour}.
\begin{claim}
\label{restriction}
Let $\mathcal{S}$ be the set of ordinals $\sigma < \lambda_i^{++}$ such that:
\begin{enumerate}
\item The restriction of $\mathcal U$ to $\widetilde{V}_\sigma$, i.e.,
\[
\mathcal{U} \cap \widetilde{V}_\sigma = \langle  \mathcal{U}(\alpha, \beta) \cap \widetilde{V}_\sigma:  \alpha \leq \lambda_i, \beta < o^{\mathcal U}(\alpha)    \rangle
\]
is a coherent sequence of measures in $\widetilde{V}_\sigma$,
\item For each $\sigma$ as above, $\alpha \leq \lambda_i$ and  $\beta < o^{\mathcal U}(\alpha)$  set
$i_{\alpha, \beta}^\sigma: V^* \to M_{\alpha, \beta}^{*, \sigma} \simeq \Ult(V^*, \mathcal{U}(\alpha, \beta) \cap \widetilde{V}_\sigma)$.
Then
the sequence
\[
\mathcal{K} \cap \widetilde{V}_\sigma= \langle   K_{\alpha, \beta} \cap \widetilde{V}_\sigma : \alpha \leq \lambda_i, \beta < o^{\mathcal U}(\alpha)         \rangle
 \]
 ia a sequence of filters such that:
 \begin{enumerate}
 \item $K_{\alpha, \beta}  \cap \widetilde{V}_\sigma$ is $\bigg(\Col(\kappa^{+4}, < i_{\alpha, \beta}(\kappa)) \times \Add(\kappa^{+3}, i_{\alpha, \beta}(\kappa)) \times  \Add(\kappa^{+4}, i_{\alpha, \beta}(\kappa)^+)\bigg)_{M^{*, \sigma}_{\alpha, \beta}}$-generic over $M^{*, \sigma}_{\alpha, \beta}$,

 \item the sequence is coherent in the sense that
 \[
 \langle K_{\alpha, \tau} \cap \widetilde{V}_\sigma : \tau < \beta  \rangle = [\bar\alpha \mapsto \langle K_{\bar\alpha, \tau} \cap \widetilde{V}_\sigma: \tau < \beta \rangle]_{\mathcal{U}(\alpha, \beta) \cap \widetilde{V}_\sigma}.
 \]
 \end{enumerate}
 \end{enumerate}
 Then $\mathcal{S}$ is a stationary subset of $\lambda_i^{++}$.
\end{claim}
Now take $\sigma \in \mathcal{S}$ such that $T \in \widetilde{V}_\sigma$. Let us define the graph $T^* \in \widetilde{V}_{\lambda_i^{++}}$
as follows:
\begin{itemize}
\item $T^*$  has the universe $\lambda_i \cup Y,$ where $Y \subseteq \lambda_i^{++}$ has size $\lambda=\lambda_i^+$, $\min(Y) > \sigma$
and $Y \cap \mathcal{S}$ is cofinal in $\sup(Y)$.

\item the pair $(\gamma, \beta)$ is an edge in $T^*$ iff  $\gamma \in Y$ and $\beta \in A_\gamma$
or symmetrically $\beta \in Y$ and $\gamma \in A_\beta.$
\end{itemize}
Since $T$ is universal, we can find and embedding $f: T^* \to T$ with $f \in \widetilde{V}_{\lambda_i^{++}}$. Take $X \subseteq \lambda_i^{++}$
of size $\leq \lambda_i$ such that
$f \restriction \lambda_i \in \widetilde{V}_\sigma[\langle A_\alpha: \alpha \in X      \rangle]$.
The same arguments as in \cite{friedman-thompson} show that there exists $\gamma \in Y \cap \mathcal{S} \setminus X$ such that $A_\gamma \notin \widetilde{V}_\sigma[\langle A_\alpha: \alpha \in X      \rangle].$
By \cite[Lemma 2.6]{friedman-thompson} we can recover $A_\gamma$ using
$f \restriction \lambda_i$ and $T$, which implies $A_\gamma \in \widetilde{V}_\sigma[\langle A_\alpha: \alpha \in X      \rangle],$
which is impossible.

\begin{remark}
In the model $W$, for every singular cardinal $\lambda,$ we have $2^{<\lambda}=\lambda$ and hence by classical results in model theory, there is a universal graph of size $\lambda$ in $W$.
\end{remark}
\begin{remark}
In the model $W$, for each uncountable regular cardinal $\lambda,$ the universality number for graphs on $\lambda$
is $\lambda^{++}$.
\end{remark}

\end{document}